\documentclass[12pt]{article}
\usepackage{amsmath}
\usepackage{amssymb}
\usepackage{amscd}
\usepackage{array}
\usepackage{amsthm}
\usepackage[flushmargin]{footmisc}
\newtheorem{thm}{Theorem}[section]

\theoremstyle{definition}

\newcommand{\Q}{\mathbb Q}

\newcommand{\Gal}{\mathrm{Gal}}

\newcommand{\Hol}{\mathrm{Hol}}

\newcommand{\Sym}{\operatorname{Sym}}

\newcommand{\eps}{\varepsilon}

\newcommand{\End}{\operatorname{End}}

\newcommand{\Aut}{\operatorname{Aut}}

\newcommand{\wk}{{\widetilde{K}}}

\title{Non-isomorphic Hopf-Galois structures with isomorphic underlying Hopf algebras.}
\author{Teresa Crespo, Anna Rio and Montserrat Vela}
\date{\today}

\begin{document}
\maketitle

\let\thefootnote\relax\footnote{2010MSC: 16T05, 12F10. \\ T. Crespo acknowledges support by grant MTM2012-33830, Spanish Science Ministry, and 2009SGR 1370; A. Rio and M. Vela acknowledge support by grant MTM2012-34611, Spanish Science Ministry, and 2009SGR 1220.}

\begin{abstract}
We give a degree 8 separable extension having two non-isomorphic
Hopf-Galois structures with isomorphic underlying Hopf algebras.
\end{abstract}

\section{Introduction}

A finite extension of fields $K/k$ is a Hopf Galois extension if
there exist a finite cocommutative $k$-Hopf algebra $H$  and a
Hopf action of $H$ on $K$, i.e a $k$-linear map $\mu: H\to
\End_k(K)$ inducing a bijection $K\otimes_k H\to\End_k(K)$. We
shall call such a pair $(H,\mu)$ a Hopf Galois structure on $K/k$.
Two Hopf Galois structures $(H_1,\mu_1)$  and $(H_2,\mu_2)$ on
$K/k$ are isomorphic if there exists a Hopf algebra isomorphism
$f:H_1 \rightarrow H_2$ such that $\mu_2 \circ f=\mu_1$.

For a Hopf Galois extension we have the following \emph{Galois
correspondence theorem}:

\begin{thm}[\cite{ChS}, Theorem 7.6]
Let $K/k$ be a finite field extension and let $(H,\mu)$ be a Hopf Galois
structure on $K/k$. For a $k$-sub-Hopf algebra $H'$ of $H$ we
define
$$
K^{H'}=\{x\in K \mid \mu(h)(x)=\eps(h)\cdot x \mbox{ for all }
h\in H'\}
$$
where $\eps$ is the counit of $H$. Then, $K^{H'}$ is a subfield
of $K$, containing $k$, and
$$
\begin{array}{rcl}
{\cal F}_H:\{H'\subseteq H \mbox{ sub-Hopf algebra}\}&\longrightarrow&\{\mbox{Fields }E\mid k\subseteq E\subseteq K\}\\
H'&\mapsto &K^{H'}
\end{array}
$$
is injective and inclusion reversing.
\end{thm}

Separable Hopf Galois extensions have been characterized by
Greither and Pareigis. Let $K/k$ be a separable extension of
degree $n$ and let $\wk/k$ be its Galois closure, $G=\Gal(\wk/k)$
and $G'=\Gal(\wk/K)$. The action of $G$ on the set
$G/G'$ of left cosets induces a group morphism

$$
 \begin{array}{rl}
\lambda: & G\rightarrow \Sym(G/G')=S_n\\
 &g\mapsto (\lambda_g: xG'\mapsto gxG')\, .
\end{array}
$$

\begin{thm}[\cite{GP} Theorem 3.1]\label{N}
With the above notations, a finite separable field extension $K/k$
is a Hopf Galois extension if and only if there exists a regular
subgroup $N$ of $S_n$ normalized by $\lambda (G)$. Moreover, there is a
one-to-one correspondence between the set of isomorphism classes
of Hopf Galois structures on $K/k$ and the set of regular
subgroups $N$ of $S_n$ normalized by $\lambda (G)$.
\end{thm}

The underlying Hopf algebra of the Hopf Galois structure corresponding to a subgroup $N$
is the sub-Hopf algebra $H=\wk[N]^G$ of $G$-fixed points in the group algebra $\wk[N]$, where $G$ acts on $\wk$ by field automorphisms and on
$N$ by conjugation inside $S_n$ via $\lambda$.

We shall use the following equivalent condition in our study of Hopf Galois structures.

\begin{thm}[\cite{Childs1} Proposition 1]
The group $\lambda (G)$ normalizes the regular subgroup $N$ of $\Sym(G/G')$ if and
only if $\lambda (G)$ is a subgroup of the holomorph $\Hol(N)=N \rtimes \Aut N$ of $N$.
\end{thm}

Within the class $\mathcal{HG}$ of separable Hopf Galois
extensions, Greither and Pareigis defined the subclass $\mathcal{AC}$ of almost classically
Galois extensions and proved that $\mathcal{AC}$ is included in the class $\mathcal{B}$
of separable Hopf Galois
extensions which may be endowed with a
Hopf Galois structure such that the Galois correspondence is
bijective (\cite{GP} \S 4). In \cite{CRV}, we proved that the
two inclusions $\mathcal{AC} \subset \mathcal{B}$ and $\mathcal{B} \subset \mathcal{HG}$ are strict.

In \cite{CRV2}, Example 2.1, we exhibited a degree 8 Hopf Galois
extension of the field $\Q$ of rational numbers which is not
almost classically Galois. In this paper, we study the Hopf Galois
structures of this extension and prove that there are two
non-isomorphic Hopf Galois structures with isomorphic underlying
Hopf algebras. In the last section, we observe that the image of
the Galois correspondence of a Hopf Galois structure does not
determine the isomorphism class of the underlying Hopf algebra.

\section{A Hopf Galois degree 8 extension}

Let us consider a Galois extension $\widetilde{K}/k$ with Galois
group $G$ isomorphic to the symmetric group $S_4$ and let $G'$ be
an order 3 subgroup of $G$. Let $K$ be the subfield of
$\widetilde{K}$ fixed under the action of $G'$. If $S_4$ is
realized over $k$ as the splitting field of a degree four
polynomial $P(X) \in k[X]$, then $K=k(\alpha,\sqrt{\delta})$,
where $\alpha$ is a root of $P(X)$ and $\delta$ its discriminant.
We showed in \cite{CRV2} Proposition 2.1 that such an extension
$K/k$ is Hopf Galois but not almost classically Galois.

\subsection{Hopf Galois structures}

According to theorem \ref{N}, in order to determine the Hopf
Galois structures of $K/k$, we look at regular subgroups $N$ of
$S_8$. Since the set of isomorphism classes of regular subgroups
of $S_8$ is in one-to-one correspondence with the set of
isomorphism classes of groups of order 8, we have five possible
groups $N$, up to isomorphism, namely the abelian groups
$C_8$, $C_2 \times C_4$ and $C_2 \times C_2 \times C_2$, the
dihedral group $D_{8}$ and the quaternion group $H_8$. If we look
for holomorphs having order divisible by 24, we are left with $C_2
\times C_2 \times C_2$ and $H_8$. But the holomorf of $H_8$ has no
transitive subgroups isomorphic to $S_4$. So the only possible
$N$'s are isomorphic to $C_2 \times C_2 \times C_2$.

Let us now determine explicitly the immersion $\lambda$ of $G$
into $Perm(G/G')$. Taking $G=\langle \tau:=(1,2,3,4),\sigma:=(1,2)\rangle$  and
$G'=\langle (2,3,4) \rangle$, a left transversal of $G'$ in $G$ is
$\{1_G,(1,2), (1,3),
(1,4), (2,3),
(1,2,3),(1,3,4), (1,4,2)
\}.$ By choosing an enumeration
of the cosets, we may take $\lambda(G)=\langle \lambda(\tau)=(1,2,3,4)(5,6,7,8),
\lambda(\sigma)=(1,2)(3,5)(4,6)(7,8)\rangle$ and then
$\lambda(G')=\langle (2,4,5)(3,8,6)\rangle$.

The conjugation class of regular subgroups of $S_8$ isomorphic to
$C_2 \times C_2 \times C_2$ has length 30. There are exactly four
subgroups $N$ in it satisfying $\lambda(G) \subset Norm_{S_8}(N)$,
namely:

$$\begin{array}{l}
N_1=\langle r_1=(1,3)(2,4)(5,7)(6,8), s_1=(1,8)(2,7)(3,6)(4,5), t=(1,7)(2,8)(3,5)(4,6)\rangle\\
N_2=\langle r_2=(1,3)(2,6)(4,8)(5,7), s_2=(1,4)(2,5)(3,8)(6,7), t=(1,7)(2,8)(3,5)(4,6)\rangle\\
N_3=\langle r_3=(1,6)(2,4)(3,8)(5,7), s_3=(1,7)(2,3)(4,8)(5,6), t_3=(1,8)(2,5)(3,6)(4,7)\rangle\\
N_4=\langle r_4=(1,3)(2,5)(4,7)(6,8),s_4=(1,7)(2,6)(3,4)(5,8),
t_4=(1,6)(2,7)(3,8)(4,5)
\rangle
\end{array}$$

We have then four Hopf Galois structures on $K/k$, up to
isomorphism.

\subsection{Galois correspondence}

By classical Galois theory, the extension $K/k$ has two strictly
intermediate fields, $k(\sqrt{\delta})$ and $k(\alpha)$.  We shall
now determine the image of the Galois correspondence for each of
the Hopf Galois structures. To this end, we compute the subgroups
of each of the corresponding $N$'s which are stable under
conjugation by $\lambda(G)$ (see the reformulation of the Galois correspondence theorem
in terms of groups in \cite{CRV} Theorem 2.3). The action of $\tau$ and $\sigma$ on
the generators of each $N_i$, $1\leq i \leq 4$, is given in the following table.

\vspace{0.5cm}
\begin{center}
\begin{tabular}{|c||c|c|c|c|c|c|c|c|c|c|c|}
\hline
& $r_1$ & $s_1$ & $t$ &$r_2$ & $s_2$ & $r_3$ & $s_3$ & $t_3$& $r_4$ & $s_4$ & $t_4$ \\
\hline \hline
$\tau$ & $r_1$ & $r_1 s_1$ & $t$ &$r_2t$ & $r_2s_2t$ & $r_3t_3$ & $r_3s_3$ & $t_3$& $r_4t_4$ & $r_4s_4$ & $t_4$ \\
\hline
$\sigma$ & $r_1s_1$ & $s_1$ & $t$ &$s_2$ & $r_2$ & $r_3$ & $r_3s_3$ & $r_3t_3$& $r_4$ & $r_4s_4$ & $r_4t_4$ \\
\hline
\end{tabular}
\end{center}

\vspace{0.5cm} The $\lambda(G)$-stable subgroups of each $N$ and
the corresponding intermediate fields of $K/k$ are given in the
following table.

\vspace{0.5cm} \begin{center}
\begin{tabular}{c|c|c}
Groups & Stable subgroups & Fixed subfields \\
\hline
$N_1$ & $\langle t\rangle, \langle r_1, s_1\rangle$ & $k(\alpha), k(\sqrt{\delta})$ \\
$N_2$ & $\langle t\rangle$ & $k(\alpha)$ \\
$N_3$ & $\langle r_3, t_3\rangle$ & $k(\sqrt{\delta})$ \\
$N_4$ & $\langle r_4, t_4\rangle$ & $k(\sqrt{\delta})$  \\
\end{tabular}
\end{center}

\subsection{A Hopf algebra with two different Hopf Galois structures}

We shall analize now if two of the underlying Hopf algebras of the four Hopf Galois structures
on $K/k$ may be isomorphic. Two Hopf algebras $H_i=\widetilde{K}[N_i]^G$ and $H_j=\widetilde{K}[N_j]^G$ are isomorphic if and
only if the groups $N_i$ and $N_j$ are $G$-isomorphic. By looking at the stable subgroups of
each of the groups $N$, we see that the only possible $G$-isomorphism is between $N_3$ and $N_4$. Now the isomorphism $\Phi:N_3  \rightarrow N_4$ defined by
$\Phi(r_3)=r_4, \Phi(s_3)=s_4, \Phi(t_3)=t_4$ is clearly a
$G$-isomorphism. It may also be seen as induced by conjugation by
$s=(1,7)(2,8)(3,5)(4,6)$  and $s$ satisfies $s g^{-1} s g \in
Cent_{S_8}(N_3)$ since $s \tau^{-1} s \tau =s \sigma^{-1} s \sigma
=Id$.  One may check that $N_3$ has no nontrivial
$G$-automorphisms, hence $\Phi$ is the unique $G$-isomorphism from
$N_3$ onto $N_4$.

The isomorphism $\Phi$ induces a Hopf algebra isomorphism between
$H_3=\widetilde{K}[N_3]^G$ and $H_4=\widetilde{K}[N_4]^G$ but
Theorem \ref{N} implies that the corresponding Hopf Galois
structures are not isomorphic. We do now a direct check of this
fact. We determine first the  Hopf algebra $H_3=\widetilde{K}[N_3]^G= \{ h \in \widetilde{K}[N_3] \, | \, \empty^g\! h=h \, , \forall g \in G \}$. From the action
of $G$ on $N_3$ determined above we obtain that an element
$h=\sum_{n \in N_3} a_n n= a_{0} Id+ a_{1} r_3 + a_{2} s_3 + a_{3}
t_3 + a_{4} r_3 s_3 + a_{5} r_3 t_3+a_6 s_3 t_3+ a_7 r_3 s_3 t_3
\in \wk[N_3]$ belongs to $H_3$ if and only if

$$\begin{array}{l} a_0 \in k, a_1 \in \wk^{\langle \sigma, \tau^2\rangle}, a_2 \in
k(\alpha), a_3=\sigma \tau(a_1), \\ a_5=\tau(a_1), a_4=\tau(a_2),
a_6=\tau^2(a_2), a_7=\tau^3(a_2).\end{array}$$

If $N$ is a regular subgroup of $Perm(G/G')$ normalized by $G$ and
$H=\wk[N]^G$ is the corresponding Hopf algebra, the Hopf action
$\mu: H \rightarrow End_k(K)$ is given by

$$(\sum a_n n) \cdot x= \sum a_n (n^{-1})(1_G)(x).$$

To make the Hopf actions $\mu_3$ and $\mu_4$ of $H_3$ and $H_4$
explicit, we first compute the preimage of $\overline{1_G}$ under the
elements of $N_3$ and $N_4$.

\vspace{0.5cm}

\begin{tabular}{c||c|c|c|c|c|c|c|c}

$n \in N_3$ & \,  $Id$ \,  & $r_3$ & $s_3$ & $t_3$ & $r_3s_3$ & $r_3t_3$ & $s_3t_3$ & $r_3s_3t_3$
\\
\hline
$(n^{-1})(\overline{1_G})$ & $\overline{1_G}$ & $\overline{(1,4,2)}$&$\overline{(2,3)}$& $\overline{(1,2,3)}$ & $\overline{(1,3)}$&$\overline{(1,3,4)}$ & $\overline{(1,4)}$&
$\overline{(1,2)}$ \end{tabular}

\vspace{0.5cm}

\begin{tabular}{c||c|c|c|c|c|c|c|c}

$n \in N_4$ &  \, $Id$ \,  & $r_4$ & $s_4$ & $t_4$ & $r_4s_4$ & $r_4t_4$ & $s_4t_4$ & $r_4s_4t_4$
\\
\hline
$(n^{-1})(\overline{1_G})$ & $\overline{1_G}$ & $\overline{(1,3,4)}$&$\overline{(2,3)}$& $\overline{(1,4,2)}$ & $\overline{(1,4)}$&$\overline{(1,2,3)}$ & $\overline{(1,2)}$&
$\overline{(1,3)}$ \end{tabular}

\vspace{0.5cm}

Let $\alpha_1=\alpha, \alpha_2, \alpha_3, \alpha_4$ be the four roots
of the polynomial $P(X)$ in $\widetilde{K}$. We consider the element $h= \alpha_1^2
s_3+ \alpha_2^2 r_3 s_3+ \alpha_3^2 s_3 t_3+\alpha_4^2 r_3 s_3 t_3
\in H_3$. We have

$$\begin{array}{rll} \mu_3(h)(\alpha_1) & = & \alpha_1^3+\alpha_2^2 \alpha_3+ \alpha_3^2 \alpha_4+\alpha_4^2
\alpha_2, \\ \mu_4(\Phi(h))(\alpha_1)&=& \mu_4(\alpha_1^2 s_4+
\alpha_2^2 r_4 s_4+ \alpha_3^2 s_4 t_4+\alpha_4^2 r_4 s_4
t_4)(\alpha_1)\\ &=& \alpha_1^3+\alpha_2^2\alpha_4+\alpha_3^2\alpha_2+\alpha_4^2\alpha_3.\end{array}$$

\noindent In order to see $\alpha_2^2 \alpha_3+ \alpha_3^2
\alpha_4+\alpha_4^2 \alpha_2 \neq
\alpha_2^2\alpha_4+\alpha_3^2\alpha_2+\alpha_4^2\alpha_3$, we
write both elements in the $k$-basis
$(\alpha_4^{i_4}\alpha_3^{i_3}\alpha_2^{i_2})_{0\leq i_4\leq 3, 0\leq i_3\leq 2,
0\leq i_2\leq 1}$ of $\wk$. Writing $P(X)=X^4-b_1X^3+b_2X^2-b_3X+b_4$, we obtain

$$\begin{array}{l}\alpha_2^2 \alpha_3+ \alpha_3^2 \alpha_4+\alpha_4^2
\alpha_2=\alpha_4^3+(\alpha_2-b_1)\alpha_4^2+(-\alpha_2\alpha_3+b_2+\alpha_3^2)
\alpha_4 +(-\alpha_3^2+b_1\alpha_3)\alpha_2-b_3, \\
\alpha_2^2\alpha_4+\alpha_3^2\alpha_2+\alpha_4^2\alpha_3=
-\alpha_4^3+(-\alpha_2+b_1)
\alpha_4^2+((b_1-\alpha_3)\alpha_2-b_2+b_1\alpha_3-
\alpha_3^2)\alpha_4+\alpha_2\alpha_3^2. \end{array}$$

\section{Final remarks}

Let $K/k$ be a separable extension of degree $n$ and let $G$ be
the Galois group of the Galois closure of $K/k$. Clearly, the fact
that two regular subgroups of $S_n$ are $G$-isomorphic implies
that the corresponding Hopf Galois structures have the same Galois
correspondence image. The converse is not true. In \cite{CRV}
Theorem 3.4, we give a family of Hopf Galois extensions having two
Hopf Galois structures of cyclic and Frobenius type, respectively,
with the same Galois correspondence image.

Let us now consider a Galois extension $K/k$ with Galois group a
Hamiltonian group $G$, i.e. $G$ is a non-abelian group such that
all its subgroups are normal subgroups. All Hamiltonian groups are
of the form  $H_8\times B \times D$, where $H_8$ is the quaternion
group, $B$ is the direct sum of some number of copies of the
cyclic group $C_2$, and D is a periodic abelian group with all
elements of odd order (see \cite{H} \S 12.5). Let $n$ be the order
of $G$. We may consider two Hopf Galois structures on $K/k$
associated to two subgroups of $S_n$ isomorphic to $G$, the group
$\rho(G)$, where $\rho$ is induced by the action of $G$ on itself
by right translation, and the group $\lambda(G)$. For the first
one, the Hopf algebra is the group algebra $k[G]$ and the Hopf
action is the linear extension of the action of the Galois group
$G$ by $k$-automorphisms of $K$. For the second one, the image of
the Galois correspondence is the set of intermediate fields $E$
such that $E/k$ is Galois (\cite{GP}, Theorem 5.3), which in the
case of Hamiltonian groups is the whole sublattice. We obtain then
two Hopf Galois structures with same Galois correspondence image
associated to two isomorphic but not $G$-isomorphic regular
subgroups of $S_n$.

\vspace{1cm} \footnotesize \noindent Teresa Crespo, Departament
d'\`Algebra i Geometria, Universitat de Barcelona, Gran Via de les
Corts Catalanes 585, E-08007 Barcelona, Spain, e-mail:
teresa.crespo@ub.edu

\vspace{0.3cm}
\noindent Anna Rio, Departament de Matem\`{a}tica Aplicada II, Universitat Polit\`{e}cnica de Catalunya, C/Jordi Girona, 1-3- Edifici Omega, E-08034 Barcelona, Spain, e-mail: ana.rio@upc.edu

\vspace{0.3cm} \noindent Montserrat Vela, Departament de
Matem\`{a}tica Aplicada II, Universitat Polit\`{e}cnica de
Catalunya, C/Jordi Girona, 1-3- Edifici Omega, E-08034 Barcelona,
Spain, e-mail: montse.vela@upc.edu
\end{document}